\documentclass[11pt]{article}
\usepackage{amsfonts}
\usepackage{amsfonts}
\usepackage{amsfonts}
\usepackage{mathrsfs}
\usepackage{bm}
\usepackage{cite}
\usepackage[colorlinks,linkcolor=blue,citecolor=blue]{hyperref}
\usepackage{amssymb,amsmath} 
 \allowdisplaybreaks

\catcode`!=11
\let\!int\int \def\int{\displaystyle\!int}
\let\!lim\lim \def\lim{\displaystyle\!lim}
\let\!sum\sum \def\sum{\displaystyle\!sum}
\let\!sup\sup \def\sup{\displaystyle\!sup}
\let\!inf\inf \def\inf{\displaystyle\!inf}
\let\!cap\cap \def\cap{\displaystyle\!cap}
\let\!max\max \def\max{\displaystyle\!max}
\let\!min\min \def\min{\displaystyle\!min}
\let\!frac\frac \def\frac{\displaystyle\!frac}
\catcode`!=12

\let\oldsection\section
\renewcommand\section{\setcounter{equation}{0}\oldsection}

\allowdisplaybreaks
\def\pf{\it{Proof.}\rm\quad}

\def\N{\mathbb{N}}

\newtheorem{defn}{Definition}[section]
\newtheorem{thm}{Theorem}[section]
\newtheorem{lem}[thm]{Lemma}
\newtheorem{cor}[thm]{Corollary}

\begin{document}
\title {\bf Some evaluation of infinite series involving trigonometric and hyperbolic functions}
\author{
{Ce Xu\thanks{Corresponding author. Email: xuce1242063253@163.com}}\\[1mm]
\small School of Mathematical Sciences, Xiamen University\\
\small Xiamen
361005, P.R. China}

\date{}
\maketitle \noindent{\bf Abstract} In this paper, by using the residue theorem and asymptotic formulas of trigonometric and hyperbolic functions at the poles, we establish many relations
involving two or more infinite series of trigonometric and hyperbolic trigonometric functions. In particular, we evaluate in closed form certain classes of infinite series containing hyperbolic
trigonometric functions, which are related to Gamma functions and $\pi$. Finally, some interesting new consequences and illustrative examples are considered.
\\[2mm]
\\[2mm]
\tableofcontents
\section{Introduction}

Infinite series involving the trigonometric and hyperbolic trigonometric functions have attracted the attention of many authors. For example, Ramanujan evaluate many infinite series of
 hyperbolic functions in his notebooks \cite{R2012} and lost notebook \cite{R1988}. We cite two example giving evaluation of series involving $\cosh(z)$. First, in Chapter 17 in his
 second notebook \cite{R2012,B1991}, Ramanujan asserted that
\[\sum\limits_{n = 1}^\infty  {\frac{{{{\left( {2n + 1} \right)}^2}}}
{{\cosh \left( {\left( {2n + 1} \right)\pi /2} \right)}}}  = \frac{{{\pi ^{3/2}}}}
{{2\sqrt 2 {\Gamma ^6}\left( {\frac{3}
{4}} \right)}}.\]
In fact, it is shown in \cite{B1991} that one can also evaluate in closed form the more general sum
\[\sum\limits_{n = 1}^\infty  {\frac{{{{\left( {2n + 1} \right)}^{2k}}}}
{{\cosh \left( {\left( {2n + 1} \right)\alpha /2} \right)}}} \]
where $k$ is a positive integer.

Second, the evaluation
\[\sum\limits_{n = 1}^\infty  {\frac{{{{\left( {2n + 1} \right)}^2}}}
{{{{\cosh }^2}\left( {\left( {2n + 1} \right)\pi /2} \right)}}}  = \frac{{{\pi ^2}}}
{{12{\Gamma ^8}\left( {\frac{3}
{4}} \right)}}\]
aries in Ramanujan's formulas for the power series coefficients of the reciprocals, or, more generally, quotients, of certain Eisenstein series \cite{BB2002}.
Berndt's papers \cite{B1977,B1978,B2004,B2016} and paper \cite{BS2017} with A. Straub,
and his books \cite{B1985,B1989,B1991} and his second, fourth book with Andrews \cite{AB2009,AB2013} contain many such results well as numerous references. Here $\Gamma \left(z\right)$ denotes
 the Gamma function, when ${\mathop{\Re}\nolimits} \left( z \right) > 0$ then
\[\Gamma \left( z \right) := \int\limits_0^\infty  {{e^{ - t}}{t^{z - 1}}dt}.\]

The motivation for this paper arises from three infinite series of Berndt's recent paper \cite{B2016}
\begin{align*}
\sum\limits_{n = 1}^\infty  {\frac{{{n^{4k - 1}}{{\left( { - 1} \right)}^{n - 1}}}}
{{\sinh \left( {\pi n} \right)}}} ,\ \sum\limits_{n = 1}^\infty  {\frac{{{n^{2k}}{{\left( { - 1} \right)}^{n - 1}}}}
{{\cosh \left( {\pi n} \right)}}}\quad {\rm and}\quad \sum\limits_{n = 1}^\infty  {\frac{{{{\left( {2n - 1} \right)}^{4k - 3}}{{\left( { - 1} \right)}^{n - 1}}}}
{{\cosh \left( {\frac{{2n - 1}}
{2}\pi } \right)}}}  \quad (k\in\N)
\end{align*}
and the methods of paper of Flajolet and Salvy \cite{FS1998}. The purposes of this paper are to establish some relations involving two or more infinite series containing hyperbolic and trigonometric
 functions by using the residue theorem. Then we use apply it to obtain a family of identities relating hyperbolic functions to Gamma functions and $\pi$.

\begin{defn} Let $\left| \theta  \right| < \pi $, we define the following eight functions
\begin{align*}
  {f_1}\left( {z,\theta } \right): = \frac{\pi }
{{\sin \left( {\pi z} \right)}}\frac{{{\pi ^2}\cosh \left( {\theta z} \right)}}
{{{{\sinh }^2}\left( {\pi z} \right)}}, \hfill \\
  {f_2}\left( {z,\theta } \right): = \frac{\pi }
{{\cos \left( {\pi z} \right)}}\frac{{{\pi ^2}\sinh \left( {\theta z} \right)}}
{{{{\cosh }^2}\left( {\pi z} \right)}}, \hfill \\
  {f_3}\left( {z,\theta } \right): = \frac{\pi }
{{\sin \left( {\pi z} \right)}}\frac{{{\pi ^2}\cosh \left( {\theta z} \right)}}
{{{{\cosh }^2}\left( {\pi z} \right)}}, \hfill \\
  {f_4}\left( {z,\theta } \right): = \frac{\pi }
{{\cos \left( {\pi z} \right)}}\frac{{{\pi ^2}\sinh \left( {\theta z} \right)}}
{{{{\sinh }^2}\left( {\pi z} \right)}}, \hfill \\
  {g_1}\left( {z,\theta } \right): = \pi \cot \left( {\pi z} \right)\frac{{{\pi ^2}\cosh \left( {\theta z} \right)}}
{{{z^2}{{\sinh }^2}\left( {\pi z} \right)}}, \hfill \\
  {g_2}\left( {z,\theta } \right): = \pi \tan \left( {\pi z} \right)\frac{{{\pi ^2}\cosh \left( {\theta z} \right)}}
{{{z^2}{{\cosh }^2}\left( {\pi z} \right)}}, \hfill \\
  {g_3}\left( {z,\theta } \right): = \pi \cot \left( {\pi z} \right)\frac{{{\pi ^2}\cosh \left( {\theta z} \right)}}
{{{z^2}{{\cosh }^2}\left( {\pi z} \right)}}, \hfill \\
  {g_4}\left( {z,\theta } \right): = \pi \tan \left( {\pi z} \right)\frac{{{\pi ^2}\cosh \left( {\theta z} \right)}}
{{{z^2}{{\sinh }^2}\left( {\pi z} \right)}}.
\end{align*}
\end{defn}
In this paper, we will consider the above eight functions and use the residue theorem to prove our main results. Here $\theta\in(-\pi,\pi)$. The functions $f_i(z,\theta)$ and $g_i(z,\theta)$ ($i=1,2,3,4$) are meromorphic in the entire complex plane.
The basic techniques of this paper, beyond the Cauchy-Lindel$\ddot{\rm o}$f contour integrals of Lemma \ref{lem1} and the asymptotic formulas of trigonometric and hyperbolic functions at the poles of Lemma \ref{lem2}, have been worked out in an experimental manner using Mathematica or Maple.

The following lemmas will be useful in the development of the main theorems.
\begin{lem}(\cite{FS1998})\label{lem1}
Let $\xi \left( s \right)$ be a kernel function and let $r(s)$ be a function which is $O(s^{-2})$ at infinity. Then
\begin{align}\label{e1}
\sum\limits_{\alpha  \in O} {{\mathop{\rm Re}\nolimits} s{{\left( {r\left( s \right)\xi \left( s \right)} \right)}_{s = \alpha }}}  + \sum\limits_{\beta  \in S} {{\mathop{\rm Re}\nolimits} s{{\left( {r\left( s \right)\xi \left( s \right)} \right)}_{s = \beta }}}  = 0.
\end{align}
where $S$ is the set of poles of $r(s)$ and $O$ is the set of poles of $\xi \left( s \right)$ that are not poles $r(s)$ . Here ${\mathop{\rm Re}\nolimits} s{\left( {r\left( s \right)} \right)_{s = \alpha }} $ denotes the residue of $r(s)$ at $s= \alpha$. The kernel function $\xi \left( s \right)$ is meromorphic in the whole complex plane and satisfies $\xi \left( s \right)=o(s)$ over an infinite collection of circles $\left| s \right| = {\rho _k}$ with ${\rho _k} \to \infty . $
\end{lem}
\pf It suffices to apply the residue theorem to
\[\frac{1}
{{2\pi i}}\int\limits_{(\infty)}  {r\left( z \right)\xi \left( z \right)} dz\]
where $\oint\limits_{\left( \infty  \right)}$ denotes integration along large circles, that is, the limit of integrals $\oint\limits_{\left| z \right| = p_k} $.\hfill$\square$

It is clear that the formula (\ref{e1}) is also true if $r\left( z \right)\xi \left( z \right)=o(z^{-\alpha})\ (\alpha>1)$ over an infinite collection of circles $\left| z \right| = {\rho _k}$ with ${\rho _k} \to \infty$.

\begin{lem}(\cite{X2017})\label{lem2} Let $n$ be a integer, then the following asymptotic formulas hold:
\begin{align}\label{e2}
\pi \cot \left( {\pi z} \right)\mathop  = \limits^{z \to n} \frac{1}{{z - n}} - 2\sum\limits_{k = 1}^\infty  {\zeta \left( {2k} \right){{\left( {z - n} \right)}^{2k - 1}}} ,
\end{align}
\begin{align}\label{e3}
\frac{\pi }
{{\sin \left( {\pi z} \right)}}\mathop  = \limits^{z \to n} {\left( { - 1} \right)^n}\left( {\frac{1}
{{z - n}} + 2\sum\limits_{k = 1}^\infty  {\bar \zeta \left( {2k} \right){{\left( {z - n} \right)}^{2k - 1}}} } \right),
\end{align}
\begin{align}\label{e4}
\pi \coth \left( {\pi z} \right)\mathop  = \limits^{z \to ni} \frac{1}
{{z - ni}} - 2\sum\limits_{k = 1}^\infty  {{{\left( { - 1} \right)}^k}\zeta \left( {2k} \right){{\left( {z - ni} \right)}^{2k - 1}}} ,
\end{align}
\begin{align}\label{e5}
\frac{\pi }
{{\sinh \left( {\pi z} \right)}}\mathop  = \limits^{z \to ni} {\left( { - 1} \right)^n}\left( {\frac{1}
{{z - ni}} + 2\sum\limits_{k = 1}^\infty  {{{\left( { - 1} \right)}^k}\bar \zeta \left( {2k} \right){{\left( {z - ni} \right)}^{2k - 1}}} } \right),
\end{align}
\begin{align}\label{e6}
 \pi \tan \left( {\pi z} \right)\mathop  = \limits^{z \to \frac{{2n - 1}}
{2}}  - \frac{1}
{{z - \frac{{2n - 1}}
{2}}} + 2\sum\limits_{k = 1}^\infty  {\zeta \left( {2k} \right){{\left( {z - \frac{{2n - 1}}
{2}} \right)}^{2k - 1}}} ,
\end{align}
\begin{align}\label{e7}
\pi \tanh \left( {\pi z} \right)\mathop  = \limits^{z \to \frac{{2n - 1}}
{2}i} \frac{1}
{{z - \frac{{2n - 1}}
{2}i}} - 2\sum\limits_{k = 1}^\infty  {{{\left( { - 1} \right)}^k}\zeta \left( {2k} \right){{\left( {z - \frac{{2n - 1}}
{2}i} \right)}^{2k - 1}}} ,
\end{align}
\begin{align}\label{e8}
\frac{\pi }
{{\cos \left( {\pi z} \right)}}\mathop  = \limits^{z \to \frac{{2n - 1}}
{2}} {\left( { - 1} \right)^n}\left\{ {\frac{1}
{{z - \frac{{2n - 1}}
{2}}} + 2\sum\limits_{k = 1}^\infty  {\bar \zeta \left( {2k} \right){{\left( {z - \frac{{2n - 1}}
{2}} \right)}^{2k - 1}}} } \right\},
\end{align}
\begin{align}\label{e9}
&   \frac{\pi }
{{\cosh \left( {\pi z} \right)}}\mathop  = \limits^{z \to \frac{{2n - 1}}
{2}i} {\left( { - 1} \right)^n}i\left\{ {\frac{1}
{{z - \frac{{2n - 1}}
{2}i}} + 2\sum\limits_{k = 1}^\infty  {{{\left( { - 1} \right)}^k}\bar \zeta \left( {2k} \right){{\left( {z - \frac{{2n - 1}}
{2}i} \right)}^{2k - 1}}} } \right\},
\end{align}
where $\zeta\left(s \right)$ and ${\bar \zeta}\left(s \right)$ denote the Riemann zeta function and alternating Riemann zeta function, respectively, which are defined by \cite{A2000}
$$\zeta(s):=\sum\limits_{n = 1}^\infty {\frac {1}{n^{s}}},\Re(s)>1,$$
\[\bar \zeta \left( s \right) := \sum\limits_{n = 1}^\infty  {\frac{{{{\left( { - 1} \right)}^{n - 1}}}}{{{n^s}}}}=(1-2^{1-s})\zeta(s) ,\;{\mathop{\Re}\nolimits} \left( s \right) \ge 1.\]
\end{lem}
When $s=2k\ (k\in \N)$ is an even, Euler proved the famous formula
\begin{align}\label{e10}
\zeta \left( {2k} \right) = \frac{{{{\left( { - 1} \right)}^{k - 1}}{B_{2k}}}}
{{2\left( {2k} \right)!}}{\left( {2\pi } \right)^{2k}},
\end{align}
where $B_{2k}$ is Bernoulli number.
The formula not only provides an elegant formula for evaluating $\zeta(2k)$, but it also tells us of the arithmetical nature of $\zeta(2k)$. In contrast, we know very little about
 $\zeta(2k+1)$. One of the major achievements in number theory in the past half-century are R. Ap$\acute{\rm e}$ry's proof that $\zeta(3)$ is irrational \cite{Ap1981}, but for $k\geq 2$,
 the arithmetical nature of $\zeta(2k+1)$ remains open. Further results for $k\geq 2$, please see Zudilin's paper \cite{Z2001,2Z2001}.

Ramanujan made a lot of nice and elegant discoveries in his short life, and one of them that has attracted the attention of several mathematicians over the years is his
intriguing formula for $\zeta(2k+1)$ (see \cite{S1,R2012,R1988})
\begin{align*}
&{\left( {4\beta } \right)^{ - \left( {k - 1} \right)}}\left\{ {\frac{1}
{2}\zeta \left( {2k - 1} \right) + \sum\limits_{n = 1}^\infty  {\frac{1}
{{{n^{2k - 1}}\left( {{e^{2n\alpha }} - 1} \right)}}} } \right\}\nonumber\\
& \quad- {\left( { - 4\alpha } \right)^{ - \left( {k - 1} \right)}}\left\{ {\frac{1}
{2}\zeta \left( {2k - 1} \right) + \sum\limits_{n = 1}^\infty  {\frac{1}
{{{n^{2k - 1}}\left( {{e^{2n\beta }} - 1} \right)}}} } \right\}\nonumber\\
& = \sum\limits_{j = 0}^k {{{\left( { - 1} \right)}^{j - 1}}\frac{{{B_{2j}}{B_{2k - 2j}}}}
{{\left( {2j} \right){\text{!}}\left( {2k - 2j} \right)!}}{\alpha ^j}{\beta ^{k - j}}},
\end{align*}
To be sure, Ramanujan¡¯s formula does not possess the elegance of (\ref{e10}), nor does it provide any arithmetical information. But,
it is indeed a remarkable formula.

The remainder of this paper is organized as follows.

In the second section we use residue theorem with the help of asymptotic formulas of trigonometric and hyperbolic functions to prove eight identities (\ref{2.1})-(\ref{2.8}). Then upon
differentiating both members of (\ref{2.1})-(\ref{2.4}), we can get four corollaries. Finally, we give some specific cases.

In the last section third section we get some known consequences and examples. Then using the results of subsection \ref{sec2.5}, we obtain some interesting new illustrative
examples.
\begin{table}[htbp]\centering
 \begin{tabular}{lll}
  \hline
  kernel function $\xi \left( z\right)$ &  base function $r\left(z\right)$ & combined function $r\left(z\right)\xi \left( z\right)$ \\
 \hline
$\frac{\pi }
{{\sin \left( {\pi z} \right)}}$ & $\frac{{{\pi ^2}\cosh \left( {\theta z} \right)}}
{{{{\sinh }^2}\left( {\pi z} \right)}}$
 & ${f_1}\left( {z,\theta } \right) = \frac{\pi }
{{\sin \left( {\pi z} \right)}}\frac{{{\pi ^2}\cosh \left( {\theta z} \right)}}
{{{{\sinh }^2}\left( {\pi z} \right)}}$
 \\
$\frac{\pi }
{{\cos \left( {\pi z} \right)}}$
 &$\frac{{{\pi ^2}\sinh \left( {\theta z} \right)}}
{{{{\cosh }^2}\left( {\pi z} \right)}}$
 &${f_2}\left( {z,\theta } \right) = \frac{\pi }
{{\cos \left( {\pi z} \right)}}\frac{{{\pi ^2}\sinh \left( {\theta z} \right)}}
{{{{\cosh }^2}\left( {\pi z} \right)}}$
\\
$\frac{\pi }
{{\sin \left( {\pi z} \right)}}$
& $\frac{{{\pi ^2}\cosh \left( {\theta z} \right)}}
{{{{\cosh }^2}\left( {\pi z} \right)}}$
& ${f_3}\left( {z,\theta } \right) = \frac{\pi }
{{\sin \left( {\pi z} \right)}}\frac{{{\pi ^2}\cosh \left( {\theta z} \right)}}
{{{{\cosh }^2}\left( {\pi z} \right)}}$
\\
$\frac{\pi }
{{\cos \left( {\pi z} \right)}}$
&$\frac{{{\pi ^2}\sinh \left( {\theta z} \right)}}
{{{{\sinh }^2}\left( {\pi z} \right)}}$
&${f_4}\left( {z,\theta } \right) = \frac{\pi }
{{\cos \left( {\pi z} \right)}}\frac{{{\pi ^2}\sinh \left( {\theta z} \right)}}
{{{{\sinh }^2}\left( {\pi z} \right)}}$
\\
$\pi \cot \left( {\pi z} \right)$
&$\frac{{{\pi ^2}\cosh \left( {\theta z} \right)}}
{{{z^2}{{\sinh }^2}\left( {\pi z} \right)}}$
&${g_1}\left( {z,\theta } \right) = \pi \cot \left( {\pi z} \right)\frac{{{\pi ^2}\cosh \left( {\theta z} \right)}}
{{{z^2}{{\sinh }^2}\left( {\pi z} \right)}}$
\\
$\pi \tan \left( {\pi z} \right)$
&$\frac{{{\pi ^2}\cosh \left( {\theta z} \right)}}
{{{z^2}{{\cosh }^2}\left( {\pi z} \right)}}$
&${g_2}\left( {z,\theta } \right) = \pi \tan \left( {\pi z} \right)\frac{{{\pi ^2}\cosh \left( {\theta z} \right)}}
{{{z^2}{{\cosh }^2}\left( {\pi z} \right)}}$\\
$\pi \cot \left( {\pi z} \right)$
&$\frac{{{\pi ^2}\cosh \left( {\theta z} \right)}}
{{{z^2}{{\cosh }^2}\left( {\pi z} \right)}}$
&${g_3}\left( {z,\theta } \right) = \pi \cot \left( {\pi z} \right)\frac{{{\pi ^2}\cosh \left( {\theta z} \right)}}
{{{z^2}{{\cosh }^2}\left( {\pi z} \right)}}$
\\
$\pi \tan \left( {\pi z} \right)$
&$\frac{{{\pi ^2}\cosh \left( {\theta z} \right)}}
{{{z^2}{{\sinh }^2}\left( {\pi z} \right)}}$
&${g_4}\left( {z,\theta } \right) = \pi \tan \left( {\pi z} \right)\frac{{{\pi ^2}\cosh \left( {\theta z} \right)}}
{{{z^2}{{\sinh }^2}\left( {\pi z} \right)}}$\\
  \hline
 \end{tabular}
 \begin{center}
  \textbf{\footnotesize\bf TABLE 2.1.}\ \ kernels and base functions
  \end{center}
\end{table}\\

\section{Main results and proofs}
In this section we state our main results. The main results of this section are the following two theorems. In the context of this paper, the results of Theorem \ref{thm2.1}-\ref{thm2.2}
 can be proved by applying the following kernels and base functions (see Table 2.1)
\subsection{Eight identities}

\begin{thm}\label{thm2.1} For any $\left| \theta  \right| < \pi $, then
\begin{align}\label{2.1}
{\pi ^2}\sum\limits_{n = 1}^\infty  {\frac{{\cosh \left( {n\theta } \right)}}
{{{{\sinh }^2}\left( {\pi n} \right)}}{{\left( { - 1} \right)}^n}}  + \theta \pi \sum\limits_{n = 1}^\infty  {\frac{{\sin \left( {n\theta } \right)}}
{{\sinh \left( {\pi n} \right)}}}  + {\pi ^2}\sum\limits_{n = 1}^\infty  {\frac{{\cos \left( {n\theta } \right)\cosh \left( {\pi n} \right)}}
{{{{\sinh }^2}\left( {\pi n} \right)}}}  + \frac{{{\theta ^2}}}
{4} - \frac{1}
{2}\zeta \left( 2 \right) = 0,
\end{align}

\begin{align}\label{2.2}
{\pi ^2}\sum\limits_{n = 1}^\infty  {\frac{{\sinh \left( {\frac{{2n - 1}}
{2}\theta } \right)}}
{{{{\cosh }^2}\left( {\frac{{2n - 1}}
{2}\pi } \right)}}{{\left( { - 1} \right)}^n}} & - \theta \pi \sum\limits_{n = 1}^\infty  {\frac{{\cos \left( {\frac{{2n - 1}}
{2}\theta } \right)}}
{{\cosh \left( {\frac{{2n - 1}}
{2}\pi } \right)}}}\nonumber \\
& + {\pi ^2}\sum\limits_{n = 1}^\infty  {\frac{{\sin \left( {\frac{{2n - 1}}
{2}\theta } \right)\sinh \left( {\frac{{2n - 1}}
{2}\pi } \right)}}
{{{{\cosh }^2}\left( {\frac{{2n - 1}}
{2}\pi } \right)}}}  = 0,
\end{align}

\begin{align}\label{2.3}
{\pi ^2}\sum\limits_{n = 1}^\infty  {\frac{{\cosh \left( {n\theta } \right)}}
{{{{\cosh }^2}\left( {\pi n} \right)}}{{\left( { - 1} \right)}^n}} & - \theta \pi \sum\limits_{n = 1}^\infty  {\frac{{\sin \left( {\frac{{2n - 1}}
{2}\theta } \right)}}
{{\sinh \left( {\frac{{2n - 1}}
{2}\pi } \right)}}} \nonumber\\
& - {\pi ^2}\sum\limits_{n = 1}^\infty  {\frac{{\cos \left( {\frac{{2n - 1}}
{2}\theta } \right)\cosh \left( {\frac{{2n - 1}}
{2}\pi } \right)}}
{{{{\sinh }^2}\left( {\frac{{2n - 1}}
{2}\pi } \right)}}}  + {\pi ^2} = 0,
\end{align}

\begin{align}\label{2.4}
{\pi ^2}\sum\limits_{n = 1}^\infty  {\frac{{\sinh \left( {\frac{{2n - 1}}
{2}\theta } \right)}}
{{{{\sinh }^2}\left( {\frac{{2n - 1}}
{2}\pi } \right)}}{{\left( { - 1} \right)}^n}}  + \theta \pi \sum\limits_{n = 1}^\infty  {\frac{{\cos \left( {n\theta } \right)}}
{{\cosh \left( {\pi n} \right)}}}  - {\pi ^2}\sum\limits_{n = 1}^\infty  {\frac{{\sin \left( {n\theta } \right)\sinh \left( {\pi n} \right)}}
{{{{\cosh }^2}\left( {\pi n} \right)}}}  + \pi \theta  = 0.
\end{align}
\end{thm}

\begin{thm}\label{thm2.2} For any $\left| \theta  \right| < \pi $, then
\begin{align}\label{2.5}
&{\pi ^2}\sum\limits_{n = 1}^\infty  {\frac{{\cosh \left( {n\theta } \right)}}
{{{n^2}{{\sinh }^2}\left( {\pi n} \right)}}}  - {\pi ^2}\sum\limits_{n = 1}^\infty  {\frac{{\cos \left( {n\theta } \right)}}
{{{n^2}{{\sinh }^2}\left( {\pi n} \right)}}}  - \theta \pi \sum\limits_{n = 1}^\infty  {\frac{{\coth \left( {\pi n} \right)\sin \left( {n\theta } \right)}}
{{{n^2}}}}\nonumber \\
&\quad- 2\pi \sum\limits_{n = 1}^\infty  {\frac{{\coth \left( {\pi n} \right)\cos \left( {n\theta } \right)}}
{{{n^3}}}}  + \frac{{{\theta ^4}}}
{{48}} - {\theta ^2}\zeta \left( 2 \right) + 7\zeta \left( 4 \right) = 0,
\end{align}

\begin{align}\label{2.6}
&{\pi ^2}\sum\limits_{n = 1}^\infty  {\frac{{\cosh \left( {\frac{{2n - 1}}
{2}\theta } \right)}}
{{{{\left( {2n - 1} \right)}^2}{{\cosh }^2}\left( {\frac{{2n - 1}}
{2}\pi } \right)}}}  - {\pi ^2}\sum\limits_{n = 1}^\infty  {\frac{{\cos \left( {\frac{{2n - 1}}
{2}\theta } \right)}}
{{{{\left( {2n - 1} \right)}^2}{{\cosh }^2}\left( {\frac{{2n - 1}}
{2}\pi } \right)}}} \nonumber\\
& + \theta \pi \sum\limits_{n = 1}^\infty  {\frac{{\tanh \left( {\frac{{2n - 1}}
{2}\pi } \right)\sin \left( {\frac{{2n - 1}}
{2}\theta } \right)}}
{{{{\left( {2n - 1} \right)}^2}}}}  + 4\pi \sum\limits_{n = 1}^\infty  {\frac{{\tanh \left( {\frac{{2n - 1}}
{2}\pi } \right)\cos \left( {\frac{{2n - 1}}
{2}\theta } \right)}}
{{{{\left( {2n - 1} \right)}^3}}}} \nonumber \\&\quad - \frac{{{\pi ^4}}}
{8} = 0,
\end{align}

\begin{align}\label{2.7}
&{\pi ^2}\sum\limits_{n = 1}^\infty  {\frac{{\cosh \left( {n\theta } \right)}}
{{{n^2}{{\cosh }^2}\left( {\pi n} \right)}}}  + 4{\pi ^2}\sum\limits_{n = 1}^\infty  {\frac{{\cos \left( {\frac{{2n - 1}}
{2}\theta } \right)}}
{{{{\left( {2n - 1} \right)}^2}{{\sinh }^2}\left( {\frac{{2n - 1}}
{2}\pi } \right)}}} \nonumber\\
& + 4\theta \pi \sum\limits_{n = 1}^\infty  {\frac{{\coth \left( {\frac{{2n - 1}}
{2}\pi } \right)\sin \left( {\frac{{2n - 1}}
{2}\theta } \right)}}
{{{{\left( {2n - 1} \right)}^2}}}}  + 16\pi \sum\limits_{n = 1}^\infty  {\frac{{\coth \left( {\frac{{2n - 1}}
{2}\pi } \right)\cos \left( {\frac{{2n - 1}}
{2}\theta } \right)}}
{{{{\left( {2n - 1} \right)}^3}}}}\nonumber\\& \quad + \left( {\frac{{{\theta ^2}}}
{4}- 4\zeta \left( 2 \right)} \right){\pi ^2} = 0,
\end{align}

\begin{align}\label{2.8}
&{\pi ^2}\sum\limits_{n = 1}^\infty  {\frac{{\cos \left( {n\theta } \right)}}
{{{n^2}{{\cosh }^2}\left( {\pi n} \right)}}}  + 4{\pi ^2}\sum\limits_{n = 1}^\infty  {\frac{{\cosh \left( {\frac{{2n - 1}}
{2}\theta } \right)}}
{{{{\left( {2n - 1} \right)}^2}{{\sinh }^2}\left( {\frac{{2n - 1}}
{2}\pi } \right)}}}  - \theta \pi \sum\limits_{n = 1}^\infty  {\frac{{\tanh \left( {\pi n} \right)\sin \left( {n\theta } \right)}}
{{{n^2}}}} \nonumber\\
& \quad- 2\pi \sum\limits_{n = 1}^\infty  {\frac{{\tanh \left( {\pi n} \right)\cos \left( {n\theta } \right)}}
{{{n^3}}}}  - \frac{{{\pi ^2}{\theta ^2}}}
{4} = 0.
\end{align}

\end{thm}
\subsection{Proofs of formulas (\ref{2.1})-(\ref{2.4})}
Now we are ready to prove Theorem \ref{thm2.1}. The proof is based on the functions $f_i(z,\theta)\ (i=1,2,3,4)$ and the usual residue computation. The functions $f_1(z,\theta),\ f_2(z,\theta),\ f_3(z,\theta)$ and $f_4(z,\theta)$ are meromorphic in the entire complex plane with some simple poles, see Table 2.2 (Here $n$ is positive integer).
\begin{table}[htbp]\centering
 \begin{tabular}{lll}
  \hline
  poles of kernel function  & poles of base function  & function  \\
 \hline
 $0,\ \pm n $& $0,\ \pm ni$ & $f_1(z,\theta)$ \\
$\pm \frac {2n-1}{2}$ & $\pm \frac {2n-1}{2}i$ &$f_2(z,\theta)$\\
$0,\ \pm n $ & $\pm \frac {2n-1}{2}i$ &$f_3(z,\theta)$\\
$\pm \frac {2n-1}{2}$ &$0,\ \pm ni$ &$f_4(z,\theta)$\\
  \hline
 \end{tabular}
 \begin{center}
  \textbf{\footnotesize\bf TABLE 2.2.}\ \ poles of kernels and base functions
  \end{center}
\end{table}

First, we prove the first identity (\ref{2.1}). The only singularities are poles at the integers and $0,\ \pm ni\ (n\in \N)$. At a integer $\pm n$ the pole is simple and the residue is
\[{\text{Res}}\left[ {{f_1}\left( {z,\theta } \right),z =  \pm n} \right] = {\left( { - 1} \right)^n}\frac{{{\pi ^2}\cosh \left( {n\theta } \right)}}
{{{{\sinh }^2}\left( {\pi n} \right)}}.\]
Then from (\ref{e5}), we can find that
\begin{align}\label{2.9}
\frac{{{\pi ^2}}}
{{{{\sinh }^2}\left( {\pi z} \right)}} = \frac{1}
{{{{\left( {z - ni} \right)}^2}}} - 4\bar \zeta \left( 2 \right) + o\left( {{{\left( {z - ni} \right)}^2}} \right),\quad n\rightarrow ni.
\end{align}
Hence, at a imaginary number $\pm ni$, the pole has order two and the residue is
$${\text{Res}}\left[ {{f_1}\left( {z,\theta } \right),z =  \pm ni} \right] = \frac{{\theta \pi \sin \left( {n\theta } \right)}}
{{\sinh \left( {\pi n} \right)}} + \frac{{{\pi ^2}\cos \left( {n\theta } \right)\cosh \left( {\pi n} \right)}}
{{{{\sinh }^2}\left( {\pi n} \right)}}.$$
Furthermore, applying (\ref{e3}) and (\ref{2.9}), we deduce that
\[{f_1}\left( {z,\theta } \right) = \cosh \left( {\theta z} \right)\left\{ {\frac{1}
{{{z^3}}} - 2\frac{{\bar \zeta \left( 2 \right)}}
{z} + o\left( z \right)} \right\},\quad z\rightarrow 0.\]
Therefore, the residue of the pole of order three at 0 is found to be
\[{\text{Res}}\left[ {{f_1}\left( {z,\theta } \right),z = 0} \right]{\text{ = }}\frac{{{\theta ^2}}}
{2} - \zeta \left( 2 \right).\]
By Lemma \ref{lem1}, summing these three contributions yields the desired result (\ref{2.1}).

Similarly, the function $f_2(z,\theta)$ only have poles at the $\pm(2n-1)/2$ and $\pm(2n-1)i/2$, $f_3(z,\theta)$ only have poles at the integer and $\pm(2n-1)i/2$, $f_4(z,\theta)$ only have poles at the $\pm(2n-1)/2$ and $0, \pm ni$.
By a simple residue computations with the help of asymptotic formulas (\ref{e2})-(\ref{e9}), we easily obtain
\begin{align*}
&{\text{Res}}\left[ {{f_2}\left( {z,\theta } \right),z =  \pm \frac{{2n - 1}}
{2}} \right] = {\left( { - 1} \right)^n}\frac{{{\pi ^2}\sinh \left( {\frac{{2n - 1}}
{2}\theta } \right)}}
{{{{\cosh }^2}\left( {\frac{{2n - 1}}
{2}\pi } \right)}},\\
&{\text{Res}}\left[ {{f_2}\left( {z,\theta } \right),z =  \pm \frac{{2n - 1}}
{2}i} \right] =  - \frac{{\theta \pi \cos \left( {\frac{{2n - 1}}
{2}\theta } \right)}}
{{\cosh \left( {\frac{{2n - 1}}
{2}\pi } \right)}} + \frac{{{\pi ^2}\sin \left( {\frac{{2n - 1}}
{2}\theta } \right)\sinh \left( {\frac{{2n - 1}}
{2}\pi } \right)}}
{{{{\cosh }^2}\left( {\frac{{2n - 1}}
{2}\pi } \right)}},\\
& {\text{Res}}\left[ {{f_3}\left( {z,\theta } \right),z = 0} \right] = {\pi ^2},  \\
& {\text{Res}}\left[ {{f_3}\left( {z,\theta } \right),z =  \pm n} \right] = {\left( { - 1} \right)^n}\frac{{{\pi ^2}\cosh \left( {n\theta } \right)}}
{{{{\cosh }^2}\left( {\pi n} \right)}},  \\
&{\text{Res}}\left[ {{f_3}\left( {z,\theta } \right),z =  \pm \frac{{2n - 1}}
{2}i} \right] =  - \frac{{\theta \pi \sin \left( {\frac{{2n - 1}}
{2}\theta } \right)}}
{{\sinh \left( {\frac{{2n - 1}}
{2}\pi } \right)}} - \frac{{{\pi ^2}\cos \left( {\frac{{2n - 1}}
{2}\theta } \right)\cosh \left( {\frac{{2n - 1}}
{2}\pi } \right)}}
{{{{\sinh }^2}\left( {\frac{{2n - 1}}
{2}\pi } \right)}},\\
&{\text{Res}}\left[ {{f_4}\left( {z,\theta } \right),z = 0} \right] = \theta \pi , \\
& {\text{Res}}\left[ {{f_4}\left( {z,\theta } \right),z =  \pm \frac{{2n - 1}}
{2}} \right] = {\left( { - 1} \right)^n}\frac{{{\pi ^2}\sinh \left( {\frac{{2n - 1}}
{2}\theta } \right)}}
{{{{\sinh }^2}\left( {\frac{{2n - 1}}
{2}\pi } \right)}}, \\
& {\text{Res}}\left[ {{f_4}\left( {z,\theta } \right),z =  \pm ni} \right] = \frac{{\theta \pi \cos \left( {n\theta } \right)}}
{{\cosh \left( {\pi n} \right)}} - \frac{{{\pi ^2}\sin \left( {n\theta } \right)\sinh \left( {\pi n} \right)}}
{{{{\cosh }^2}\left( {\pi n} \right)}}.
\end{align*}
Hence, by Lemma \ref{lem1}, by a direct calculation we may easily deduce the formulas (\ref{2.2})-(\ref{2.4}). This completes the proof of Theorem \ref{thm2.1}. \hfill$\square$

\subsection{Proofs of formulas (\ref{2.5})-(\ref{2.8})}
The proofs of formulas (\ref{2.5})-(\ref{2.8}) are similar to the proofs of formulas (\ref{2.1})-(\ref{2.4}). These results from applying the residue theorem to the functions $g_i\ (i=1,2,3,4)$.
First, we can get the poles of functions $g_i\ (i=1,2,3,4)$, see Table 2.3.
\begin{table}[htbp]\centering
 \begin{tabular}{lll}
  \hline
  poles of kernel function  & poles of base function  & function  \\
 \hline
 $0,\ \pm n $& $0,\ \pm ni$ & $g_1(z,\theta)$ \\
$\pm \frac {2n-1}{2}$ & 0,\ $\pm \frac {2n-1}{2}i$ &$g_2(z,\theta)$\\
$0,\ \pm n $ & $0,\ \pm \frac {2n-1}{2}i$ &$g_3(z,\theta)$\\
$\pm \frac {2n-1}{2}$ &$0,\ \pm ni$ &$g_4(z,\theta)$\\
  \hline
 \end{tabular}
 \begin{center}
  \textbf{\footnotesize\bf TABLE 2.3.}\ \ poles of kernels and base functions
  \end{center}
\end{table}\\
Then by a similar argument as in the proof of formula (\ref{2.1}) with help of asymptotic formulas (\ref{e2})-(\ref{e9}), we deduce that
\begin{align*}
 & {\text{Res}}\left[ {{g_1}\left( {z,\theta } \right),z = 0} \right] = \frac{{{\theta ^4}}}
{{24}} - 2{\theta ^2}\zeta \left( 2 \right) + 14\zeta \left( 4 \right), \hfill \\
 & {\text{Res}}\left[ {{g_1}\left( {z,\theta } \right),z =  \pm n} \right] = \frac{{{\pi ^2}\cosh \left( {n\theta } \right)}}
{{{n^2}{{\sinh }^2}\left( {\pi n} \right)}}, \hfill \\
 & {\text{Res}}\left[ {{g_1}\left( {z,\theta } \right),z =  \pm ni} \right] =  - \frac{{{\pi ^2}\cos \left( {n\theta } \right)}}
{{{n^2}{{\sinh }^2}\left( {\pi n} \right)}} - \theta \pi \frac{{\coth \left( {\pi n} \right)\sin \left( {n\theta } \right)}}
{{{n^2}}} - 2\pi \frac{{\coth \left( {\pi n} \right)\cos \left( {n\theta } \right)}}
{{{n^3}}}, \hfill \\
 & {\text{Res}}\left[ {{g_2}\left( {z,\theta } \right),z = 0} \right] = {\pi ^4}, \hfill \\
 & {\text{Res}}\left[ {{g_2}\left( {z,\theta } \right),z =  \pm \frac{{2n - 1}}
{2}} \right] =  - 4\frac{{{\pi ^2}\cosh \left( {\frac{{2n - 1}}
{2}\theta } \right)}}
{{{{\left( {2n - 1} \right)}^2}{{\cosh }^2}\left( {\frac{{2n - 1}}
{2}\pi } \right)}}, \hfill \\
&{\text{Res}}\left[ {{g_2}\left( {z,\theta } \right),z =  \pm \frac{{2n - 1}}
{2}i} \right] = 4\frac{{{\pi ^2}\cos \left( {\frac{{2n - 1}}
{2}\theta } \right)}}
{{{{\left( {2n - 1} \right)}^2}{{\cosh }^2}\left( {\frac{{2n - 1}}
{2}\pi } \right)}} - 4\theta \pi \frac{{\tanh \left( {\frac{{2n - 1}}
{2}\pi } \right)\sin \left( {\frac{{2n - 1}}
{2}\theta } \right)}}
{{{{\left( {2n - 1} \right)}^2}}}\\
&\quad\quad\quad\quad\quad\quad\quad\quad\quad\quad\quad\quad\quad\quad\quad - 16\pi \frac{{\tanh \left( {\frac{{2n - 1}}
{2}\pi } \right)\cos \left( {\frac{{2n - 1}}
{2}\theta } \right)}}
{{{{\left( {2n - 1} \right)}^3}}},\\
&{\text{Res}}\left[ {{g_3}\left( {z,\theta } \right),z = 0} \right] = \left( {\frac{{{\theta ^2}}}
{2} - 8\zeta \left( 2 \right)} \right){\pi ^2}, \hfill \\
 & {\text{Res}}\left[ {{g_3}\left( {z,\theta } \right),z =  \pm n} \right] = \frac{{{\pi ^2}\cosh \left( {n\theta } \right)}}
{{{n^2}{{\cosh }^2}\left( {\pi n} \right)}}, \hfill \\
 & {\text{Res}}\left[ {{g_3}\left( {z,\theta } \right),z =  \pm \frac{{2n - 1}}
{2}i} \right] = 4\frac{{{\pi ^2}\cos \left( {\frac{{2n - 1}}
{2}\theta } \right)}}
{{{{\left( {2n - 1} \right)}^2}{{\sinh }^2}\left( {\frac{{2n - 1}}
{2}\pi } \right)}} + 4\theta \pi \frac{{\coth \left( {\frac{{2n - 1}}
{2}\pi } \right)\sin \left( {\frac{{2n - 1}}
{2}\theta } \right)}}
{{{{\left( {2n - 1} \right)}^2}}} \hfill \\
& \quad\quad\quad\quad\quad\quad\quad\quad\quad\quad\quad\quad\quad\quad\quad + 16\pi \frac{{\coth \left( {\frac{{2n - 1}}
{2}\pi } \right)\cos \left( {\frac{{2n - 1}}
{2}\theta } \right)}}
{{{{\left( {2n - 1} \right)}^3}}},\\
&{\text{Res}}\left[ {{g_4}\left( {z,\theta } \right),z = 0} \right] = \frac{{{\pi ^2}{\theta ^2}}}
{2}, \hfill \\
&  {\text{Res}}\left[ {{g_4}\left( {z,\theta } \right),z =  \pm \frac{{2n - 1}}
{2}} \right] =  - 4\frac{{{\pi ^2}\cosh \left( {\frac{{2n - 1}}
{2}\theta } \right)}}
{{{{\left( {2n - 1} \right)}^2}{{\sinh }^2}\left( {\frac{{2n - 1}}
{2}\pi } \right)}}, \hfill \\
 & {\text{Res}}\left[ {{g_4}\left( {z,\theta } \right),z =  \pm ni} \right] =  - \frac{{{\pi ^2}\cos \left( {n\theta } \right)}}
{{{n^2}{{\cosh }^2}\left( {\pi n} \right)}}{\text{ + }}\theta \pi \frac{{\tanh \left( {\pi n} \right)\sin \left( {n\theta } \right)}}
{{{n^2}}}{\text{ + }}2\pi \frac{{\tanh \left( {\pi n} \right)\cos \left( {n\theta } \right)}}
{{{n^3}}}{\text{.}}
\end{align*}
Thus, by applying Lemma \ref{lem1} and combining related equations, we obtain the formulas (\ref{2.5})-(\ref{2.8}). The proof of Theorem \ref{thm2.2} is finished.  \hfill$\square$

\subsection{Further corollaries}
\begin{cor}\label{cor2.3} For positive integer $k$ and $\left| \theta  \right| < \pi $, then
\begin{align}\label{2.10}
&{\pi ^2}\sum\limits_{n = 1}^\infty  {\frac{{{n^{2k}}\cosh \left( {n\theta } \right)}}
{{{{\sinh }^2}\left( {\pi n} \right)}}{{\left( { - 1} \right)}^n}}  + {\left( { - 1} \right)^k}\theta \pi \sum\limits_{n = 1}^\infty  {\frac{{{n^{2k}}\sin \left( {n\theta } \right)}}
{{\sinh \left( {\pi n} \right)}}}  + 2k{\left( { - 1} \right)^{k - 1}}\pi \sum\limits_{n = 1}^\infty  {\frac{{{n^{2k - 1}}\cos \left( {n\theta } \right)}}
{{\sinh \left( {\pi n} \right)}}}\nonumber \\
&\quad\quad\quad+ {\pi ^2}{\left( { - 1} \right)^k}\sum\limits_{n = 1}^\infty  {\frac{{{n^{2k}}\cos \left( {n\theta } \right)\cosh \left( {\pi n} \right)}}
{{{{\sinh }^2}\left( {\pi n} \right)}}}  + \left\{ {\begin{array}{*{20}{c}}
   {\frac{1}
{2},\;k = 1}  \\
   {0,\;k > 1}  \\
 \end{array} } \right. = 0,
\end{align}
\begin{align}\label{2.11}
&{\pi ^2}\sum\limits_{n = 1}^\infty  {\frac{{{n^{2k - 1}}\sinh \left( {n\theta } \right)}}
{{{{\sinh }^2}\left( {\pi n} \right)}}{{\left( { - 1} \right)}^n}}  + {\left( { - 1} \right)^{k - 1}}\theta \pi \sum\limits_{n = 1}^\infty  {\frac{{{n^{2k - 1}}\cos \left( {n\theta } \right)}}
{{\sinh \left( {\pi n} \right)}}}  + \left( {2k - 1} \right){\left( { - 1} \right)^{k - 1}}\pi \sum\limits_{n = 1}^\infty  {\frac{{{n^{2k - 2}}\sin \left( {n\theta } \right)}}
{{\sinh \left( {\pi n} \right)}}}\nonumber \\
&\quad\quad+ {\pi ^2}{\left( { - 1} \right)^k}\sum\limits_{n = 1}^\infty  {\frac{{{n^{2k - 1}}\sin \left( {n\theta } \right)\cosh \left( {\pi n} \right)}}
{{{{\sinh }^2}\left( {\pi n} \right)}}}  + \left\{ {\begin{array}{*{20}{c}}
   {\frac{\theta }
{2},\;k = 1}  \\
   {0,\;k > 1}  \\
 \end{array} } \right. = 0.
\end{align}
\end{cor}
\pf Differentiating (\ref{2.1}) $2k$ and $2k-1$ times with respect to $\theta$, respectively, and applying the following well-known identities
\begin{align*}
 & {\left( {\theta \sin \left( {n\theta } \right)} \right)^{\left( {2k} \right)}} = {\left( { - 1} \right)^k}{n^{2k}}\theta \sin \left( {n\theta } \right)
   + 2k{\left( { - 1} \right)^{k - 1}}{n^{2k - 1}}\cos \left( {n\theta } \right),  \\
 & {\left( {\theta \sin \left( {n\theta } \right)} \right)^{\left( {2k - 1} \right)}} = {\left( { - 1} \right)^{k - 1}}{n^{2k - 1}}\theta \cos \left( {n\theta } \right)
  + \left( {2k - 1} \right){\left( { - 1} \right)^{k - 1}}{n^{2k - 2}}\sin \left( {n\theta } \right),
\end{align*}
we obtain formulas (\ref{2.10}) and (\ref{2.11}).\hfill$\square$
\begin{cor}\label{cor2.4} For positive integer $k$ and $\left| \theta  \right| < \pi $, then
\begin{align}\label{2.12}
&\pi \sum\limits_{n = 1}^\infty  {\frac{{{{\left( {2n - 1} \right)}^{2k}}\sinh \left( {\frac{{2n - 1}}
{2}\theta } \right)}}
{{{{\cosh }^2}\left( {\frac{{2n - 1}}
{2}\pi } \right)}}{{\left( { - 1} \right)}^n}}  - \theta {\left( { - 1} \right)^k}\sum\limits_{n = 1}^\infty  {\frac{{{{\left( {2n - 1} \right)}^{2k}}\cos \left( {\frac{{2n - 1}}
{2}\theta } \right)}}
{{\cosh \left( {\frac{{2n - 1}}
{2}\pi } \right)}}}\nonumber \\
& \quad- 4k{\left( { - 1} \right)^k}\sum\limits_{n = 1}^\infty  {\frac{{{{\left( {2n - 1} \right)}^{2k - 1}}\sin \left( {\frac{{2n - 1}}
{2}\theta } \right)}}
{{\cosh \left( {\frac{{2n - 1}}
{2}\pi } \right)}}} \nonumber\\ &\quad + \pi {\left( { - 1} \right)^k}\sum\limits_{n = 1}^\infty  {\frac{{{{\left( {2n - 1} \right)}^{2k}}\sin \left( {\frac{{2n - 1}}
{2}\theta } \right)\sinh \left( {\frac{{2n - 1}}
{2}\pi } \right)}}
{{{{\cosh }^2}\left( {\frac{{2n - 1}}
{2}\pi } \right)}}}  = 0.
\end{align}
\end{cor}
\pf Differentiating (\ref{2.2}) $2k$ times with respect to $\theta$ and using the equality
\begin{align*}
{\left( {\theta \cos \left( {\frac{{2n - 1}}
{2}\theta } \right)} \right)^{\left( {2k} \right)}} =& {\left( { - 1} \right)^k}{\left( {\frac{{2n - 1}}
{2}} \right)^{2k}}\theta \cos \left( {\frac{{2n - 1}}
{2}\theta } \right) \\&+ 2k{\left( { - 1} \right)^k}{\left( {\frac{{2n - 1}}
{2}} \right)^{2k - 1}}\sin \left( {\frac{{2n - 1}}
{2}\theta } \right),
\end{align*}
we deduce the desired result.\hfill$\square$
\begin{cor}\label{cor2.5} For positive integer $k$ and $\left| \theta  \right| < \pi $, then
\begin{align}\label{2.13}
&{2^{2k}}\pi \sum\limits_{n = 1}^\infty  {\frac{{{n^{2k}}\cosh \left( {n\theta } \right)}}
{{{{\cosh }^2}\left( {\pi n} \right)}}{{\left( { - 1} \right)}^n}}  - {\left( { - 1} \right)^k}\theta \sum\limits_{n = 1}^\infty  {\frac{{{{\left( {2n - 1} \right)}^{2k}}\sin \left( {\frac{{2n - 1}}
{2}\theta } \right)}}
{{\sinh \left( {\frac{{2n - 1}}
{2}\pi } \right)}}}\nonumber \\
&\quad + 4k{\left( { - 1} \right)^k}\sum\limits_{n = 1}^\infty  {\frac{{{{\left( {2n - 1} \right)}^{2k - 1}}\cos \left( {\frac{{2n - 1}}
{2}\theta } \right)}}
{{\sinh \left( {\frac{{2n - 1}}
{2}\pi } \right)}}} \nonumber\\
&\quad  - {\left( { - 1} \right)^k}\pi \sum\limits_{n = 1}^\infty  {\frac{{{{\left( {2n - 1} \right)}^{2k}}\cos \left( {\frac{{2n - 1}}
{2}\theta } \right)\cosh \left( {\frac{{2n - 1}}
{2}\pi } \right)}}
{{{{\sinh }^2}\left( {\frac{{2n - 1}}
{2}\pi } \right)}}}  = 0.
\end{align}
\end{cor}
\pf Similarly as in the proofs of Corollary \ref{cor2.3} and \ref{cor2.4}, we differentiate (\ref{2.3}) $2k$ times with respect to $\theta$, and with the help of formula
\[{\left( {\theta \sin \left( {\frac{{2n - 1}}
{2}\theta } \right)} \right)^{\left( {2k} \right)}} = {\left( { - 1} \right)^k}{\left( {\frac{{2n - 1}}
{2}} \right)^{2k}}\theta \sin \left( {\frac{{2n - 1}}
{2}\theta } \right) + 2k{\left( { - 1} \right)^{k - 1}}{\left( {\frac{{2n - 1}}
{2}} \right)^{2k - 1}}\cos \left( {\frac{{2n - 1}}
{2}\theta } \right),\]
by a simple calculation, we may easily obtain (\ref{2.13}).\hfill$\square$

\begin{cor}\label{cor2.6} For positive integer $k$ and $\left| \theta  \right| < \pi $, then
\begin{align}\label{2.14}
&\pi \sum\limits_{n = 1}^\infty  {\frac{{{{\left( {2n - 1} \right)}^{2k - 1}}\cosh \left( {\frac{{2n - 1}}
{2}\theta } \right)}}
{{{{\sinh }^2}\left( {\frac{{2n - 1}}
{2}\pi } \right)}}{{\left( { - 1} \right)}^n}}  + {2^{2k - 1}}{\left( { - 1} \right)^k}\theta \sum\limits_{n = 1}^\infty  {\frac{{{n^{2k - 1}}\sin \left( {n\theta } \right)}}
{{\cosh \left( {\pi n} \right)}}} \nonumber\\
&\quad + {2^{2k - 1}}\left( {2k - 1} \right){\left( { - 1} \right)^{k - 1}}\sum\limits_{n = 1}^\infty  {\frac{{{n^{2k - 2}}\cos \left( {n\theta } \right)}}
{{\cosh \left( {\pi n} \right)}}} \nonumber\\&\quad - {2^{2k - 1}}{\left( { - 1} \right)^{k - 1}}\pi \sum\limits_{n = 1}^\infty  {\frac{{{n^{2k - 1}}\cos \left( {n\theta } \right)\sinh \left( {\pi n} \right)}}
{{{{\cosh }^2}\left( {\pi n} \right)}}}  + \left\{ {\begin{array}{*{20}{c}}
   {{1},\;k = 1}  \\
   {0,\;k > 1}  \\
 \end{array} } \right. = 0.
\end{align}
\end{cor}
\pf Differentiating (\ref{2.4}) $2k-1$ times with respect to $\theta$ and using
\[{\left( {\theta \cos \left( {n\theta } \right)} \right)^{\left( {2k - 1} \right)}} = {\left( { - 1} \right)^k}{n^{2k - 1}}\theta \sin \left( {n\theta } \right) + \left( {2k - 1} \right){\left( { - 1} \right)^{k - 1}}{n^{2k - 2}}\cos \left( {n\theta } \right),\]
by a direct calculation, the result is (\ref{2.14}).\hfill$\square$
\subsection{Some simple examples}\label{sec2.5}

\begin{align}
&\sum\limits_{n = 1}^\infty  {\frac{{\pi n{{\left( { - 1} \right)}^n}}}
{{\sinh \left( {\pi n} \right)}}}  = \frac{1}
{4},\\
&\sum\limits_{n = 1}^\infty  {\frac{{\pi {n^{4k + 1}}{{\left( { - 1} \right)}^n}}}
{{\sinh \left( {\pi n} \right)}}}  = 0,\quad (k\in \N)\\
&\sum\limits_{n = 1}^\infty  {\frac{{{\pi ^2}\cosh \left( {\pi n} \right)}}
{{{{\sinh }^2}\left( {\pi n} \right)}}{{\left( { - 1} \right)}^n} =  - \frac{1}
{2}\zeta \left( 2 \right)} ,\\
&\sum\limits_{n = 1}^\infty  {\frac{{{{\left( {2n - 1} \right)}^{4k - 1}}}}
{{\cosh \left( {\frac{{2n - 1}}
{2}\pi } \right)}}{{\left( { - 1} \right)}^n}}  = 0,\quad (k\in \N)\\
&\sum\limits_{n = 1}^\infty  {\frac{{{\pi ^2}{{\left( { - 1} \right)}^n}}}
{{{{\sinh }^2}\left( {\pi n} \right)}}}  + \sum\limits_{n = 1}^\infty  {\frac{{{\pi ^2}\cosh \left( {\pi n} \right)}}
{{{{\sinh }^2}\left( {\pi n} \right)}}} - \frac{1}
{2}\zeta \left( 2 \right) = 0,\\
&\sum\limits_{n = 1}^\infty  {\frac{{{{\left( { - 1} \right)}^n}}}
{{\cosh \left( {\pi n} \right)}}}  + \sum\limits_{n = 1}^\infty  {\frac{{{{\left( { - 1} \right)}^n}}}
{{\sinh \left( {\frac{{2n - 1}}
{2}\pi } \right)}}}  + \frac{1}
{2} = 0,\\
&\sum\limits_{n = 1}^\infty  {\frac{{{{\left( { - 1} \right)}^n}}}
{{{{\cosh }^2}\left( {\pi n} \right)}}}  - \sum\limits_{n = 1}^\infty  {\frac{{\cosh \left( {\frac{{2n - 1}}
{2}\pi } \right)}}
{{{{\sinh }^2}\left( {\frac{{2n - 1}}
{2}\pi } \right)}}}  + \frac{1}
{2} = 0,\\
&\sum\limits_{n = 1}^\infty  {\frac{{{\pi ^2}{n^2}{{\left( { - 1} \right)}^n}}}
{{{{\sinh }^2}\left( {\pi n} \right)}}}  + 2\sum\limits_{n = 1}^\infty  {\frac{{\pi n}}
{{\sinh \left( {\pi n} \right)}}}  - \sum\limits_{n = 1}^\infty  {\frac{{{\pi ^2}{n^2}\cosh \left( {\pi n} \right)}}
{{{{\sinh }^2}\left( {\pi n} \right)}}}  + \frac{1}
{2} = 0,\\
&\sum\limits_{n = 1}^\infty  {\frac{{{\pi ^2}{n^{4k}}\cosh \left( {\pi n} \right)}}
{{{{\sinh }^2}\left( {\pi n} \right)}}{{\left( { - 1} \right)}^n}}  - 2k\sum\limits_{n = 1}^\infty  {\frac{{\pi {n^{4k - 1}}{{\left( { - 1} \right)}^n}}}
{{\sinh \left( {\pi n} \right)}}}  = 0,\quad (k\in \N)\\
&{2^{2k}}\sum\limits_{n = 1}^\infty  {\frac{{{n^{2k}}{{\left( { - 1} \right)}^n}}}
{{\cosh \left( {\pi n} \right)}}}  + {\left( { - 1} \right)^k}\sum\limits_{n = 1}^\infty  {\frac{{{{\left( {2n - 1} \right)}^{2k}}}}
{{\sinh \left( {\frac{{2n - 1}}
{2}\pi } \right)}}} {\left( { - 1} \right)^n} = 0,\quad(k\in\N)\\
&4\pi \sum\limits_{n = 1}^\infty  {\frac{{{n^2}{{\left( { - 1} \right)}^n}}}
{{{{\cosh }^2}\left( {\pi n} \right)}}}  - 4\sum\limits_{n = 1}^\infty  {\frac{{2n - 1}}
{{\sinh \left( {\frac{{2n - 1}}
{2}\pi } \right)}}}  + \pi \sum\limits_{n = 1}^\infty  {\frac{{{{\left( {2n - 1} \right)}^2}\cosh \left( {\frac{{2n - 1}}
{2}\pi } \right)}}
{{{{\sinh }^2}\left( {\frac{{2n - 1}}
{2}\pi } \right)}}}  = 0,\\
&\pi \sum\limits_{n = 1}^\infty  {\frac{{{{\left( {2n - 1} \right)}^{4k - 2}}\sinh \left( {\frac{{2n - 1}}
{2}\pi } \right)}}
{{{{\cosh }^2}\left( {\frac{{2n - 1}}
{2}\pi } \right)}}{{\left( { - 1} \right)}^n}}  = 2\left( {2k - 1} \right)\sum\limits_{n = 1}^\infty  {\frac{{{{\left( {2n - 1} \right)}^{4k - 3}}}}
{{\cosh \left( {\frac{{2n - 1}}
{2}\pi } \right)}}{{\left( { - 1} \right)}^n}} ,\quad (k\in\N)\\
&\sum\limits_{n = 1}^\infty  {\frac{{\pi \cosh \left( {\pi n} \right)}}
{{{n^2}{{\sinh }^2}\left( {\pi n} \right)}}}  - \sum\limits_{n = 1}^\infty  {\frac{{\pi {{\left( { - 1} \right)}^n}}}
{{{n^2}{{\sinh }^2}\left( {\pi n} \right)}}}  - 2\sum\limits_{n = 1}^\infty  {\frac{{\coth \left( {\pi n} \right)}}
{{{n^3}}}{{\left( { - 1} \right)}^n}}  - \frac{{49}}
{{720}}{\pi ^3} = 0,\\
&\sum\limits_{n = 1}^\infty  {\frac{\pi }
{{{{\left( {2n - 1} \right)}^2}\cosh \left( {\frac{{2n - 1}}
{2}\pi } \right)}}}  + \pi \sum\limits_{n = 1}^\infty  {\frac{{\tanh \left( {\frac{{2n - 1}}
{2}\pi } \right)}}
{{{{\left( {2n - 1} \right)}^2}}}{{\left( { - 1} \right)}^{n - 1}}}  = \frac{{{\pi ^3}}}
{8},\\
&\sum\limits_{n = 1}^\infty  {\frac{\pi }
{{{n^2}\cosh \left( {\pi n} \right)}}}  + 4\pi \sum\limits_{n = 1}^\infty  {\frac{{\coth \left( {\frac{{2n - 1}}
{2}\pi } \right)}}
{{{{\left( {2n - 1} \right)}^2}}}{{\left( { - 1} \right)}^{n - 1}}}  = \frac{5}
{{12}}{\pi ^3},\\
&\sum\limits_{n = 1}^\infty  {\frac{\pi }
{{{n^2}{{\cosh }^2}\left( {\pi n} \right)}}}  + 4\sum\limits_{n = 1}^\infty  {\frac{\pi }
{{{{\left( {2n - 1} \right)}^2}{{\sinh }^2}\left( {\frac{{2n - 1}}
{2}\pi } \right)}}}  - 2\sum\limits_{n = 1}^\infty  {\frac{{\tanh \left( {\pi n} \right)}}
{{{n^3}}}}  = 0,\\
&\sum\limits_{n = 1}^\infty  {\frac{\pi }
{{{n^2}{{\cosh }^2}\left( {\pi n} \right)}}}  + 4\sum\limits_{n = 1}^\infty  {\frac{\pi }
{{{{\left( {2n - 1} \right)}^2}{{\sinh }^2}\left( {\frac{{2n - 1}}
{2}\pi } \right)}}}  + 16\sum\limits_{n = 1}^\infty  {\frac{{\coth \left( {\frac{{2n - 1}}
{2}\pi } \right)}}
{{{{\left( {2n - 1} \right)}^3}}}}  = \frac{2}
{3}{\pi ^3},\\
&\sum\limits_{n = 1}^\infty  {\frac{{\pi {{\left( { - 1} \right)}^n}}}
{{{n^2}{{\cosh }^2}\left( {\pi n} \right)}}}  + 4\sum\limits_{n = 1}^\infty  {\frac{{\pi \cosh \left( {\frac{{2n - 1}}
{2}\pi } \right)}}
{{{{\left( {2n - 1} \right)}^2}{{\sinh }^2}\left( {\frac{{2n - 1}}
{2}\pi } \right)}}}  - 2\sum\limits_{n = 1}^\infty  {\frac{{\tanh \left( {\pi n} \right)}}
{{{n^3}}}{{\left( { - 1} \right)}^n}}  = \frac{{{\pi ^3}}}
{4}.
\end{align}
These identities can be obtained from the main theorems and corollaries which are presented in the paper.
Some above results are already in the literature, e.g., the second and fourth equations were first established by Cauchy \cite{C1889}. It can also be found in Sect. 14 of Chapter 14 in
 Ramanujan's second notebook \cite{R2012,B1989} and Berndt's paper \cite{B1978,B2016}. References to several further proofs of the second and fourth equations may be found in \cite{B1989}.

\section{Some known or new illustrative examples}

\subsection{Some known results}\label{sec3.1}
In his book \cite{B1991} and paper \cite{B2016}, Berndt established many realtions between Jacobian elliptic functions and infinite series containing trigonometric and hyperbolic trigonometric
 functions. Using Berndt's method, many infinite trigonometric and hyperbolic series can be expressed in terms of Gamma functions and $\pi$. From \cite{B1991,B2016}, it is not difficult to find that
  the following
 infinite series
\begin{align*}
&\sum\limits_{n = 1}^\infty  {\frac{{{n^{4k - 1}}{{\left( { - 1} \right)}^{n - 1}}}}
{{\sinh \left( {\pi n} \right)}}} ,\ \sum\limits_{n = 1}^\infty  {\frac{{{n^{2k}}{{\left( { - 1} \right)}^{n - 1}}}}
{{\cosh \left( {\pi n} \right)}}} ,\ \sum\limits_{n = 1}^\infty  {\frac{{{n^{2k + 1}}}}
{{\sinh \left( {\pi n} \right)}}},\ \sum\limits_{n = 1}^\infty  {\frac{{{n^{2k}}}}
{{\cosh \left( {\pi n} \right)}}}\\
&\quad {\rm and}\quad \sum\limits_{n = 1}^\infty  {\frac{{{{\left( {2n - 1} \right)}^{4k - 3}}{{\left( { - 1} \right)}^{n - 1}}}}
{{\cosh \left( {\frac{{2n - 1}}
{2}\pi } \right)}}}  \quad (k\in\N)
\end{align*}
can be expressed as Gamma functions and $\pi$. We have
\begin{align}
&  \sum\limits_{n = 1}^\infty  {\frac{{{n^{4k - 1}}{{\left( { - 1} \right)}^{n - 1}}}}
{{\sinh \left( {\pi n} \right)}}}  = {c_1}\frac{{{\Gamma ^{8k}}\left( {\frac{1}
{4}} \right)}}
{{{\pi ^{6k}}}},\\
&  \sum\limits_{n = 1}^\infty  {\frac{{{n^{2k}}{{\left( { - 1} \right)}^{n - 1}}}}
{{\cosh \left( {\pi n} \right)}}}  = {c_2}\frac{{\sqrt 2 {\Gamma ^{4k + 2}}\left( {\frac{1}
{4}} \right)}}
{{{\pi ^{3k + 3/2}}}},\\
&  \sum\limits_{n = 1}^\infty  {\frac{{{n^{2k + 1}}}}
{{\sinh \left( {\pi n} \right)}}}  = {c_3}\frac{{{\Gamma ^{4k + 4}}\left( {\frac{1}
{4}} \right)}}
{{{\pi ^{3k + 3}}}},  \\
&  \sum\limits_{n = 1}^\infty  {\frac{{{n^{2k}}}}
{{\cosh \left( {\pi n} \right)}}}  = {c_4}\frac{{{\Gamma ^{4k + 2}}\left( {\frac{1}
{4}} \right)}}
{{{\pi ^{3k + 3/2}}}}, \\
& \sum\limits_{n = 1}^\infty  {\frac{{{{\left( {2n - 1} \right)}^{4k - 3}}{{\left( { - 1} \right)}^{n - 1}}}}
{{\cosh \left( {\frac{{2n - 1}}
{2}\pi } \right)}}}  = {c_5}\frac{{{\Gamma ^{4k - 2}}\left( {\frac{1}
{4}} \right)}}
{{{\pi ^{2k - 1}}{\Gamma ^{4k - 2}}\left( {\frac{3}
{4}} \right)}}
\end{align}
where the coefficients $c_1,\ c_2,\ c_3,\ c_4$ and $c_5$ are rational numbers. Using Mathematica or Maple tools, each specific coefficients are can be obtained. We can get the following examples
\begin{align}
 & \sum\limits_{n = 1}^\infty  {\frac{{\pi n}}
{{\sinh \left( {\pi n} \right)}}}  =  - \frac{1}
{4} + \frac{{{\Gamma ^4}\left( {\frac{1}
{4}} \right)}}
{{32{\pi ^2}}}, \\
 & \sum\limits_{n = 1}^\infty  {\frac{{{n^3}{{\left( { - 1} \right)}^{n - 1}}}}
{{\sinh \left( {\pi n} \right)}}}  = \frac{{{\Gamma ^8}\left( {\frac{1}
{4}} \right)}}
{{512{\pi ^6}}}, \\
 & \sum\limits_{n = 1}^\infty  {\frac{{{n^7}{{\left( { - 1} \right)}^{n - 1}}}}
{{\sinh \left( {\pi n} \right)}}}  =  - \frac{{9{\Gamma ^{16}}\left( {\frac{1}
{4}} \right)}}
{{65536{\pi ^{12}}}}, \\
 & \sum\limits_{n = 1}^\infty  {\frac{{{n^{11}}{{\left( { - 1} \right)}^{n - 1}}}}
{{\sinh \left( {\pi n} \right)}}}  = \frac{{189{\Gamma ^{24}}\left( {\frac{1}
{4}} \right)}}
{{2097152{\pi ^{18}}}},\\
&  \sum\limits_{n = 1}^\infty  {\frac{{{{\left( { - 1} \right)}^{n - 1}}}}
{{\cosh \left( {\pi n} \right)}}}  = \frac{1}
{2} - \frac{{\sqrt 2 {\Gamma ^2}\left( {\frac{1}
{4}} \right)}}
{{8{\pi ^{3/2}}}},  \\
 & \sum\limits_{n = 1}^\infty  {\frac{{{{{n^2}\left( { - 1} \right)}^{n - 1}}}}
{{\cosh \left( {\pi n} \right)}}}  = \frac{{\sqrt 2 {\Gamma ^6}\left( {\frac{1}
{4}} \right)}}
{{256{\pi ^{9/2}}}}, \\
&  \sum\limits_{n = 1}^\infty  {\frac{{{{{n^4}\left( { - 1} \right)}^{n - 1}}}}
{{\cosh \left( {\pi n} \right)}}}  = \frac{{3\sqrt 2 {\Gamma ^{10}}\left( {\frac{1}
{4}} \right)}}
{{8192{\pi ^{15/2}}}},  \\
 & \sum\limits_{n = 1}^\infty  {\frac{{{{{n^6}\left( { - 1} \right)}^{n - 1}}}}
{{\cosh \left( {\pi n} \right)}}}  =  - \frac{{27\sqrt 2 {\Gamma ^{14}}\left( {\frac{1}
{4}} \right)}}
{{262144{\pi ^{21/2}}}}, \\
&  \sum\limits_{n = 1}^\infty  {\frac{{{n^8}{{\left( { - 1} \right)}^{n - 1}}}}
{{\cosh \left( {\pi n} \right)}}}  =  - \frac{{441{\Gamma ^{18}}\left( {\frac{1}
{4}} \right)}}
{{8388608{\pi ^{27/2}}}}, \\
&  \sum\limits_{n = 1}^\infty  {\frac{1}
{{\cosh \left( {\pi n} \right)}}}  =  - \frac{1}
{2} + \frac{{{\Gamma ^2}\left( {\frac{1}
{4}} \right)}}
{{4{\pi ^{3/2}}}}, \\
&  \sum\limits_{n = 1}^\infty  {\frac{{{n^2}}}
{{\cosh \left( {\pi n} \right)}}}  = \frac{{{\Gamma ^6}\left( {\frac{1}
{4}} \right)}}
{{128{\pi ^{9/2}}}},  \\
&  \sum\limits_{n = 1}^\infty  {\frac{{{n^4}}}
{{\cosh \left( {\pi n} \right)}}}  = \frac{{9{\Gamma ^{10}}\left( {\frac{1}
{4}} \right)}}
{{4096{\pi ^{15/2}}}}, \\
&  \sum\limits_{n = 1}^\infty  {\frac{{{n^6}}}
{{\cosh \left( {\pi n} \right)}}}  =  - \frac{{153{\Gamma ^{14}}\left( {\frac{1}
{4}} \right)}}
{{131072{\pi ^{21/2}}}},  \\
&  \sum\limits_{n = 1}^\infty  {\frac{{{n^8}}}
{{\cosh \left( {\pi n} \right)}}}  = \frac{{4977{\Gamma ^{18}}\left( {\frac{1}
{4}} \right)}}
{{4194304{\pi ^{27/2}}}},\\
&  \sum\limits_{n = 1}^\infty  {\frac{{{n^3}}}
{{\sinh \left( {\pi n} \right)}}}  = \frac{{{\Gamma ^8}\left( {\frac{1}
{4}} \right)}}
{{256{\pi ^6}}}, \\
&  \sum\limits_{n = 1}^\infty  {\frac{{{n^5}}}
{{\sinh \left( {\pi n} \right)}}}  = \frac{{3{\Gamma ^{12}}\left( {\frac{1}
{4}} \right)}}
{{2048{\pi ^9}}}, \\
&  \sum\limits_{n = 1}^\infty  {\frac{{{n^7}}}
{{\sinh \left( {\pi n} \right)}}}  = \frac{{9{\Gamma ^{16}}\left( {\frac{1}
{4}} \right)}}
{{8192{\pi ^{12}}}},  \\
&  \sum\limits_{n = 1}^\infty  {\frac{{{n^9}}}
{{\sinh \left( {\pi n} \right)}}}  = \frac{{189{\Gamma ^{20}}\left( {\frac{1}
{4}} \right)}}
{{131072{\pi ^{15}}}}\\
& \sum\limits_{n = 1}^\infty  {\frac{{\left( {2n - 1} \right){{\left( { - 1} \right)}^{n - 1}}}}
{{\cosh \left( {\frac{{2n - 1}}
{2}\pi } \right)}}}  = \frac{{{\Gamma ^2}\left( {\frac{1}
{4}} \right)}}
{{8\pi {\Gamma ^2}\left( {\frac{3}
{4}} \right)}}, \\
 & \sum\limits_{n = 1}^\infty  {\frac{{{{\left( {2n - 1} \right)}^5}{{\left( { - 1} \right)}^{n - 1}}}}
{{\cosh \left( {\frac{{2n - 1}}
{2}\pi } \right)}}}  =  - \frac{{3{\Gamma ^6}\left( {\frac{1}
{4}} \right)}}
{{32{\pi ^3}{\Gamma ^6}\left( {\frac{3}
{4}} \right)}},\\
 & \sum\limits_{n = 1}^\infty  {\frac{{{{\left( {2n - 1} \right)}^9}{{\left( { - 1} \right)}^{n - 1}}}}
{{\cosh \left( {\frac{{2n - 1}}
{2}\pi } \right)}}}  = \frac{{189{\Gamma ^{10}}\left( {\frac{1}
{4}} \right)}}
{{128{\pi ^5}{\Gamma ^{10}}\left( {\frac{3}
{4}} \right)}}.
\end{align}
\subsection{Some interesting results}

From subsection \ref{sec2.5} and \ref{sec3.1}, we have the conclusion: for positive integer $k$, then
\begin{align*}
&\sum\limits_{n = 1}^\infty  {\frac{{{n^{4k}}\cosh \left( {\pi n} \right)}}
{{{{\sinh }^2}\left( {\pi n} \right)}}{{\left( { - 1} \right)}^{n - 1}}}  = {a_1}\frac{{{\Gamma ^{8k}}\left( {\frac{1}
{4}} \right)}}
{{{\pi ^{6k + 1}}}},\\
&\sum\limits_{n = 1}^\infty  {\frac{{{{\left( {2n - 1} \right)}^{2k - 2}}}}
{{\sinh \left( {\frac{{2n - 1}}
{2}\pi } \right)}}{{\left( { - 1} \right)}^{n - 1}}}  = {a_3}\frac{{\sqrt 2 {\Gamma ^{4k - 2}}\left( {\frac{1}
{4}} \right)}}
{{{\pi ^{3k - 3/2}}}},\\
&\sum\limits_{n = 1}^\infty  {\frac{{{{\left( {2n - 1} \right)}^{4k - 2}}\sinh \left( {\frac{{2n - 1}}
{2}\pi } \right)}}
{{{{\cosh }^2}\left( {\frac{{2n - 1}}
{2}\pi } \right)}}{{\left( { - 1} \right)}^{n - 1}}}  = {a_2}\frac{{{\Gamma ^{4k - 2}}\left( {\frac{1}
{4}} \right)}}
{{{\pi ^{2k}}{\Gamma ^{4k - 2}}\left( {\frac{3}
{4}} \right)}},
\end{align*}
where $a_1,a_2,a_3$ are rational numbers.  We can obtain the following interesting results
\begin{align}
 & \sum\limits_{n = 1}^\infty  {\frac{{{n^4}\cosh \left( {\pi n} \right)}}
{{{{\sinh }^2}\left( {\pi n} \right)}}{{\left( { - 1} \right)}^{n - 1}}}  = \frac{{{\Gamma ^8}\left( {\frac{1}
{4}} \right)}}
{{256{\pi ^7}}}, \hfill \\
 & \sum\limits_{n = 1}^\infty  {\frac{{{n^8}\cosh \left( {\pi n} \right)}}
{{{{\sinh }^2}\left( {\pi n} \right)}}{{\left( { - 1} \right)}^{n - 1}}}  =  - \frac{{9{\Gamma ^{16}}\left( {\frac{1}
{4}} \right)}}
{{16384{\pi ^{13}}}}, \hfill \\
 & \sum\limits_{n = 1}^\infty  {\frac{{{n^{12}}\cosh \left( {\pi n} \right)}}
{{{{\sinh }^2}\left( {\pi n} \right)}}{{\left( { - 1} \right)}^{n - 1}}}  = \frac{{567{\Gamma ^{24}}\left( {\frac{1}
{4}} \right)}}
{{1048576{\pi ^{19}}}}, \hfill \\
&\sum\limits_{n = 1}^\infty  {\frac{{{{\left( { - 1} \right)}^{n - 1}}}}
{{\sinh \left( {\frac{{2n - 1}}
{2}\pi } \right)}}}  = \frac{{\sqrt 2 {\Gamma ^2}\left( {\frac{1}
{4}} \right)}}
{{8{\pi ^{3/2}}}}, \hfill \\
 & \sum\limits_{n = 1}^\infty  {\frac{{{{\left( {2n - 1} \right)}^2}}}
{{\sinh \left( {\frac{{2n - 1}}
{2}\pi } \right)}}{{\left( { - 1} \right)}^{n - 1}}}  = \frac{{\sqrt 2 {\Gamma ^6}\left( {\frac{1}
{4}} \right)}}
{{64{\pi ^{9/2}}}}, \hfill \\
 & \sum\limits_{n = 1}^\infty  {\frac{{{{\left( {2n - 1} \right)}^4}}}
{{\sinh \left( {\frac{{2n - 1}}
{2}\pi } \right)}}{{\left( { - 1} \right)}^{n - 1}}}  =  - \frac{{3\sqrt 2 {\Gamma ^{10}}\left( {\frac{1}
{4}} \right)}}
{{512{\pi ^{15/2}}}}, \hfill \\
 & \sum\limits_{n = 1}^\infty  {\frac{{{{\left( {2n - 1} \right)}^6}}}
{{\sinh \left( {\frac{{2n - 1}}
{2}\pi } \right)}}{{\left( { - 1} \right)}^{n - 1}}}  =  - \frac{{27\sqrt 2 {\Gamma ^{10}}\left( {\frac{1}
{4}} \right)}}
{{4096{\pi ^{21/2}}}}, \hfill \\
&  \sum\limits_{n = 1}^\infty  {\frac{{{{\left( {2n - 1} \right)}^2}\sinh \left( {\frac{{2n - 1}}
{2}\pi } \right)}}
{{{{\cosh }^2}\left( {\frac{{2n - 1}}
{2}\pi } \right)}}{{\left( { - 1} \right)}^{n - 1}}}  = \frac{{{\Gamma ^2}\left( {\frac{1}
{4}} \right)}}
{{4{\pi ^2}{\Gamma ^2}\left( {\frac{3}
{4}} \right)}}, \hfill \\
 & \sum\limits_{n = 1}^\infty  {\frac{{{{\left( {2n - 1} \right)}^6}\sinh \left( {\frac{{2n - 1}}
{2}\pi } \right)}}
{{{{\cosh }^2}\left( {\frac{{2n - 1}}
{2}\pi } \right)}}{{\left( { - 1} \right)}^{n - 1}}}  =  - \frac{{9{\Gamma ^6}\left( {\frac{1}
{4}} \right)}}
{{16{\pi ^4}{\Gamma ^6}\left( {\frac{3}
{4}} \right)}}, \hfill \\
 & \sum\limits_{n = 1}^\infty  {\frac{{{{\left( {2n - 1} \right)}^{10}}\sinh \left( {\frac{{2n - 1}}
{2}\pi } \right)}}
{{{{\cosh }^2}\left( {\frac{{2n - 1}}
{2}\pi } \right)}}{{\left( { - 1} \right)}^{n - 1}}}  = \frac{{945{\Gamma ^{10}}\left( {\frac{1}
{4}} \right)}}
{{64{\pi ^6}{\Gamma ^{10}}\left( {\frac{3}
{4}} \right)}}.
\end{align}

It is possible that some of similar identities involving trigonometric and hyperbolic functions can be obtained by using the methods and techniques of the present paper. For instance,
considering the function \[F\left( {z,{\theta _1},{\theta _2}} \right): = \frac{\pi }
{{\sin \left( {\pi z} \right)}}\frac{{{\pi ^2}\cos \left( {{\theta _1}z} \right)\cos \left( {{\theta _2}z} \right)}}
{{{{\sinh }^2}\left( {\pi z} \right)}},\quad (\theta_1,\theta_2\in (-\pi,\pi))
\]
and using residue theorem, we deduce that
\begin{align}
&{\pi ^2}\sum\limits_{n = 1}^\infty  {\frac{{\cos \left( {{\theta _1}n} \right)\cos \left( {{\theta _2}n} \right)}}
{{{{\sinh }^2}\left( {\pi n} \right)}}{{\left( { - 1} \right)}^n}}  - \pi {\theta _1}\sum\limits_{n = 1}^\infty  {\frac{{\sinh \left( {{\theta _1}n} \right)\cos \left( {{\theta _2}n} \right)}}
{{\sinh \left( {\pi n} \right)}}}  + \pi {\theta _2}\sum\limits_{n = 1}^\infty  {\frac{{\cosh \left( {{\theta _1}n} \right)\sin \left( {{\theta _2}n} \right)}}
{{\sinh \left( {\pi n} \right)}}}\nonumber \\
&\quad + {\pi ^2}\sum\limits_{n = 1}^\infty  {\frac{{\cosh \left( {{\theta _1}n} \right)\cos \left( {{\theta _2}n} \right)\cosh \left( {\pi n} \right)}}
{{{{\sinh }^2}\left( {\pi n} \right)}}}  + \frac{{\theta _1^2 - \theta _2^2}}
{4} - \frac{1}
{2}\zeta \left( 2 \right) = 0.
\end{align}
\\
{\bf Acknowledgments.}  We thank the anonymous referee for suggestions which led to improvements in the exposition.

 \end{document}